\documentclass[12pt]{article}
\usepackage[english]{babel}
\usepackage[backend=bibtex]{biblatex} 
\usepackage{setspace}
\usepackage{amsmath}

\usepackage{amsthm}
\usepackage{latexsym, amssymb}
\usepackage[letterpaper,top=2cm,bottom=2cm,left=3cm,right=3cm,marginparwidth=1.75cm]{geometry}
\usepackage{graphicx}
\usepackage{algorithm}
\usepackage{algpseudocode}
\usepackage{hhline}

\usepackage[colorlinks=true, allcolors=blue]{hyperref}

\graphicspath{{./}}
\usepackage{tabularx, booktabs}

\theoremstyle{plain}
\newtheorem{thm}{Theorem}[section]
\newtheorem{lem}{Lemma}[section]
\newtheorem{prop}[thm]{Proposition}
\newtheorem{cor}{Corollary}[section]

\theoremstyle{definition}
\newtheorem{defn}{Definition}[section]

\newtheorem{exmp}{Example}[section]

\theoremstyle{remark}
\newtheorem*{rem}{Remark}

\algnewcommand\Continue{\textbf{continue}}
\algnewcommand{\True}{\textbf{true}}
\algnewcommand{\False}{\textbf{false}}

\newcommand{\CR}{\nonumber\\}

\title{Pure Simplicial and Clique Complexes with a Fixed Number of Facets}
\author{Kassahun H Betre$^1$, Yan X Zhang$^1$,  Carter Edmond$^1$\\1 San Jose State University, 1 Washington Square, Sci 148, San Jose, CA 95192}

\begin{document}
\maketitle

\begin{abstract}

We study structural and enumerative aspects of pure simplicial complexes and clique complexes. We prove a necessary and sufficient condition for any simplicial complex to be a clique complex that depends only on the list of facets. We also prove a theorem that a class of ``triangle-intersection free" pure clique complexes are uniquely determined up to isomorphism merely from the facet-adjacency matrix. Lastly, we count the number of pure simplicial complexes with a fixed number of facets and find an upper bound to the number of pure clique complexes.

\end{abstract}
\clearpage
\def\thefootnote{\fnsymbol{footnote}}
\setcounter{footnote}{0}
\tableofcontents
\newpage
\setcounter{page}{1}

\section{Introduction}
An {\em abstract simplicial complex} $(V,\Delta)$ is an ordered pair of a set $V$ (the {\em base set}) together with a collection $\Delta$ of subsets of $V$ closed under taking subsets (i.e., if $\sigma\in \Delta$ and $\tau\subset\sigma$, then $\tau\in\Delta$). Elements of $\Delta$ are the {\em faces} or {\em simplices} of the complex. The {\em dimension} of a face is its cardinality minus $1$. We consider the empty set as a face with dimension $-1$. The {\em facets} of the complex are the maximal simplices, i.e., those simplices in $\Delta$ that are not proper subsets of other simplices. For $p$ a positive integer, a simplicial complex is called $p$-{\em pure} if all its facets have the same cardinality of $p$. The number of facets of a simplicial complex is its {\em facet-order}. When it is clear from context, we may drop the ``facet" and simply say the {\em order} of a complex to refer to its facet-order. 

If $G$ is a finite simple graph, the {\em clique complex} (or  Whitney complex) of $G$ is the abstract simplicial complex with the vertex set $V(G)$ of the graph as the base set and the cliques of $G$ as the simplices. Clique complexes are the same as {\em flag} complexes: simplicial complexes $X(V,\Delta)$ satisfying the condition of ``2-determinability", i.e., the condition that if any $A\subseteq V$ has the property that if every pair of vertices $\{i,j\}\subseteq A$ are faces of $X$, then $A$ is itself a face of $X$. In other words, a simplicial complex is a flag complex if all its minimal non-faces have size $2$.

Pure simplicial complexes have wide-ranging applications in computational combinatorics, algebraic topology, and algebraic geometry. Part of the motivation for this investigation is application in  combinatorial topology and emergent quantum gravity. In quantum gravity, it has been argued that not only geometry, but also topology may be stochastic near spacetime singularities \cite{wheeler1955geons}. In such cases computation of observables requires one to compute expectation values over all geometries and topologies of a fixed spacetime dimension \cite{gibbons1977action}.  Random pure simplicial complexes provide an important model for stochastic topology and geometry. 

Nonetheless, pure complexes with a fixed number of facets are natural combinatorial objects to study in and of themselves. As an example, one may start with $q$ identical equilateral triangles and consider all shapes that can be formed by joining the pieces along edges or vertices. The set of all such shapes then forms the space of all 2-dimensional pure simplicial complexes with $q$ facets and labeled vertices. We want to know how many such objects exist. Clearly this is a very large space and one way to study aspects of it is to map the complex to its {\em facet-adjacency matrix}, which, for a complex of facet-order $q$ is a $q\times q$ symmetric matrix $Q$ whose entries $Q_{ij}$ give the size of the intersection of facets $i$ and $j$. In section $2$ we find two necessary (but not sufficient) conditions that ensure that a square symmetric integer matrix is the facet-adjacency matrix of some pure complex. In section $3$ we investigate the question of the extent to which the facet-adjacency matrix uniquely determine the complex. After first characterizing clique complexes in terms of their facets, we prove that certain types of ``triangle intersection free" clique complexes are uniquely determined up to isomorphism from their facet-adjacency matrix. We also prove that it is impossible to expect too much; specifically, we prove that information up to any finite $k$ ``levels'' (of which the facet-adjacency matrix contains up through level $k=2$) can never determine the complex with no additional information.

In Section~\ref{sec:counting} we turn our attention to counting pure simplicial and clique complexes with a fixed facet-order. We derive an explicit formula for the number of pure simplicial complexes. We also derive an upper bound for the number of pure clique complexes with a fixed number of facets. We summarize the main results and offer concluding remarks in Section~\ref{sec:conclusion}.

\section{Facet-incidence and facet-adjacency matrices}

\subsection{Labeling a Simplicial Complex; Facet-incidence Matrices}

\begin{defn}
A simplicial complex with $n$ vertices and $q$ facets is said to be {\em vertex-labeled} if the vertices are distinguished via labels from a set, say $[n]=\{1,2,\ldots,n\}$. Similarly, the complex is said to be \emph{facet-labeled} if the $q$ facets are labeled, for example, by using the set $[q]=\{1,\dots,q\}$.
\end{defn}

We can uniquely represent a vertex-labeled simplicial complex by listing its facets. The base set can then be found by taking the union of the facets. The remainder of the faces are given by the non-empty subsets of the facets. For example, the simplicial complex in Fig.~\ref{fig:cliquecomplex} can be given as the set of facets $\{\{1,2,3\},\{1,3,4\}\}$. From now on we will usually use this representation.

\begin{defn}
Consider a facet-labeled simplicial complex of facet-order  $q$. Let $n$ be the number of vertices of the complex and assume the vertices are also labeled. The {\em facet-incidence matrix} $B$ for this vertex- and facet-labeled complex is the $q\times n$ matrix where the matrix elements $B_{ij}$ are 
\begin{align}
B_{ij} = \begin{cases}
1&\text{if facet $i$ contains vertex $j$;}\\
0&\text{otherwise.}
\end{cases}
\end{align}
\end{defn}

\begin{defn}
We define a $q \times n$ matrix with $q \geq 1$ containing only $0$'s and $1$'s to be \emph{realizable} if it is the facet-incidence matrix of some facet-labeled complex.
\end{defn}
It is easy to see that such a matrix is realizable if and only if each row contains at least one $1$ and no two rows are comparable in the set-inclusion ordering (in other words, no two rows correspond to sets $A$ and $B$ such that $A \subseteq B$ or $B \subseteq A$). The first condition qualifies each row as a face and the second condition ensures that the faces can all simultaneously be maximal, and thus facets.

\begin{defn}
We define the \emph{vertex data} of a facet-labeled complex with $q$ facets $\{1, \ldots, q\}$ to be a map $c:\{0,1\}^{q} \rightarrow \{0, 1, 2, \ldots\}$ such that $c(a_1\cdots a_q)$ counts the number of vertices that belong to exactly the facets with indices $i$ such that $a_i = 1$.  
\end{defn}
By default $c(0,0,\ldots,0)=0$, and 
\begin{align}
    \sum_{a_1,\dots,a_q}c(a_1,\dots,a_q) = n,
\end{align} where $n$ is the total number of vertices in the complex. For example, for a complex with two facets, $c(0,0)=0, c(0,1)=$ the number of vertices belonging only to the second facet, $c(1,0) = $number of vertices belonging only to the first facet, and $c(1,1) = $ the number of vertices belonging only in the intersection of facets 1 and 2. 
\begin{exmp}
For the simplicial complex with facet-incidence matrix given by 
\[B=\left(\begin{matrix}
    1&1&1&0\\1&0&1&1
\end{matrix}\right),\]
$c(0,0)=0, c(0,1)=1, c(1,0)=1, c(1,1)=2$.
\end{exmp}
While it is clear that any facet-labeled complex uniquely determines its vertex data, it is not always the case that any map $c:\{0,1\}^q\rightarrow \{0,1,2,\dots\}$ corresponds to a facet-labeled complex (for a quick counterexample, consider $q=2$ and 
\[
c = \{(0,0) \rightarrow 0, (1,0) \rightarrow 0, (0,1) \rightarrow 0, (1,1) \rightarrow 1\},
\]
which would represent a complex with a single vertex belonging to both facets, but the facets must be different). We call a $c$ that can arise from an actual complex \emph{realizable}.
\begin{defn}
A map $c:\{0,1\}^q\rightarrow \{0,1,2,\dots\}$ is said to be \emph{realizable} if it is the vertex data of some simplicial complex.    
\end{defn} 
The following proposition is related to our claim above on whether or not a $q \times n$ matrix is realizable.

\begin{prop}
\label{prop:vertex-data-reconstruction}
A map $c:\{0,1\}^q \rightarrow \{0,1,2,\ldots\}$ with $q \geq 2$ is realizable if and only if for all pairs $i \neq j$, $i, j \in \{1, 2, \ldots, q\}$, there exist a $q-$tuple $(a_1, \ldots, a_q)\in\{0,1\}^q$ such that
\[
a_i = 1, a_j = 0,\; c(a_1\cdots a_q) > 0.
\]
\end{prop}
\begin{proof}
Given such a realizable map $c$, it is the vertex data of a complex, so we can construct a $q \times n$ matrix $B$ with 
\[
n = \sum_{a_1, \ldots, a_q \in \{0,1\}^q} c(a_1\cdots a_q)
\]
by creating $n$ columns of (vertical) vectors $(a_1, \ldots, a_q)$. This is just a formalization of labeling all the $n$ vertices in some order, which does not affect the feasibility of the existence of such a complex. Therefore, $c$ is realizable (as vertex data) if and only if $B$ is realizable (as a matrix). For any two facets $i \neq j$, we have facet $i$ not being a subset of facet $j$ if and only if there exists some vertex (equivalently, column) $k$ where $B_{ik} = 1$ and $B_{jk} = 0$. This condition translates to the condition seen in our statement where $c(a_1 \cdots a_q) > 0$ for $a_i = 1$ and $a_j = 0$ and for some choice of the other $a_*$. Necessitating the condition for all pairs $i$ and $j$ finishes our proof.
\end{proof}

\subsection{The Facet-adjacency Matrix}

\begin{figure}
\centering
\includegraphics[width=0.8\textwidth]{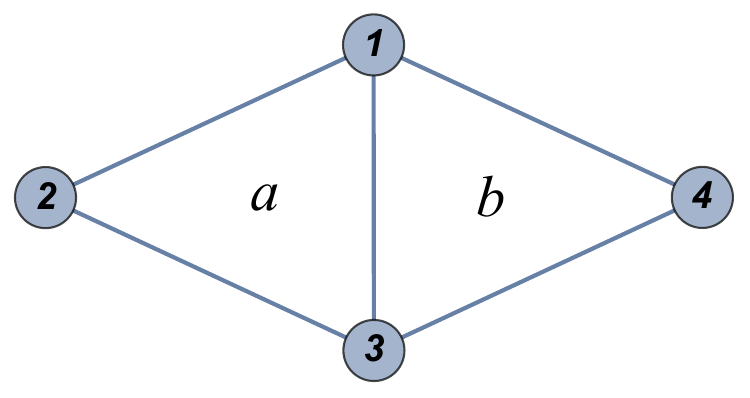}
\caption{A simplicial complex with two labeled facets and four labeled vertices.}
\label{fig:cliquecomplex}
\end{figure}

\begin{defn}
If $B$ is the facet-incidence matrix of a complex, we call the $q\times q$ square symmetric matrix $Q = BB^T$ the {\em facet-adjacency matrix}.
\end{defn} 
The value of $Q_{ij}$ is the size of the intersection of facets $i$ and $j$, so the diagonals are the facet cardinalities. All entries on row $i$ and column $j$ have values between $0$ and $Q_{ii}-1$. 
The matrices $B$ and $Q$ corresponding to the 3-pure complex in Fig. \ref{fig:cliquecomplex} are:

\begin{align}
B&=\left(\begin{matrix}
1&1&1&0\\
1&0&1&1
\end{matrix}\right),
 \quad Q=BB^T=\left(\begin{matrix} 3&2\\2&3\end{matrix} \right)
\end{align}

The facet-adjacency matrix assigns a square symmetric integer matrix to each simplicial complex. We can think of it as a map from the set of facet-labeled  complexes of facet order $q$ to $\mathbb{Z}_{q\times q}$, the set of $q\times q$ symmetric matrices with integer entries. An interesting question is the extent to which the facet-adjacency matrix characterizes a complex uniquely. It is clear that the adjacency  map is neither injective nor surjective. There are  distinct simplicial complexes with the same facet-adjacency matrix, and there are square symmetric integer matrices that are not the facet-adjacency matrices of any simplicial complex. 
\begin{defn}
    We call a $q\times q$ square symmetric matrix of integers \emph{realizable} if it is the facet-adjacency matrix of some simplicial complex.
\end{defn}

A very natural question is: \textbf{when is a matrix $Q$ realizable?} Here is one approach towards an answer:

\begin{defn}
We define a \emph{basic} matrix to be a $q \times q$ matrix such that it: 
\begin{itemize}
    \item has only $0$'s and $1$'s;
    \item is symmetric;
    \item has rank $1$.
\end{itemize}
\end{defn}

\begin{prop}
    \label{prop:fam-characterization} If $Q$ is $q \times q$  realizable, then it is a nonnegative integral linear combination of basic matrices.
\end{prop}
\begin{proof}
Let $B$ be the facet-incidence matrix of the complex and let $b_i$ denote $q \times 1$ columns of $B$, we can rewrite the facet-adjacency matrix $Q=BB^T$ as follows:

\begin{align*}
Q & = BB^T = \begin{bmatrix}b_1 & b_2 & \hdots & b_n\end{bmatrix} \begin{bmatrix}b_1^T \\ b_2^T \\ \vdots \\ b_n^T\end{bmatrix} \\
& = b_1 b_1^T + b_2 b_2^T + \cdots + b_n b_n^T.
\end{align*}

Intuitively, this means we can think about the facet-adjacency matrix $Q$ as a sum over matrices $b_ib_i^T$ corresponding to the vertices, and each vertex $i$ (belonging in $k$ facets) contributes (via $b_i b_i^T$)  with a $k \times k$ subset of (not necessarily adjacent) rows and columns containing $1$'s and $0$'s. This is exactly what it means to be a basic matrix.  
\end{proof}

\begin{cor}
        A $q\times q$ square symmetric matrix $Q$ is realizable if and only if 
    \begin{align}
        Q=\sum_{i=0}^{2^q-1}c(b_i)b_ib_i^T
    \end{align}
where $c$ is a realizable map (equivalently, the vertex data) and $b_i$ are the the $q\times q$ basic matrices.
\end{cor}
\begin{proof}
By Proposition \ref{prop:vertex-data-reconstruction}, given the vertex data of the complex we can construct the facet-incidence matrix $B$. Then, 
\begin{align}
    Q = BB^T = \sum_{i=0}^{2^q-1}c(b_i)b_ib_i^T\nonumber 
\end{align}
\end{proof}
Polytope/cone theory (see e.g. \cite{grunbaum1967convex}) gives a geometric interpretation of our situation: if a $q \times q$ matrix $Q$ is realizable, we can think of its entries as a vector in $\mathbb{Z}^{q^2}$. Proposition~\ref{prop:fam-characterization} then says that these matrices ``almost''\footnote{We say ``almost'' because we must exclude certain corner cases, such as  
$\begin{bmatrix}
   2 & 2 \\ 2 & 2 
\end{bmatrix} = 2 \begin{bmatrix}
   1 & 1 \\ 1 & 1 
\end{bmatrix}$ (which is a linear combination of basic matrices but would require two identical facets that contain the same two vertices) and $\begin{bmatrix}
   1 & 0 \\ 0 & 0 
\end{bmatrix}$ (where one of the facets must be empty).} form a cone generated by nonnegative multiples of basic matrices. The $(2^q-1)$ basic matrices, seen as vertices in $\mathbb{Z}^{q^2}$, then give a {\em vertex representation ($V$-representation)}. Equivalently, one can consider looking at the {\em half-space representation ($H$-representation)} and characterize the space by half-spaces along facets instead; these would correspond to inequalities such as the following necessary condition:
\begin{prop}  
\label{prop:inequality-constructability}
    Let $Q$ be a $q \times q$ realizable matrix. Then, for any $S\subseteq [q]$ and $i\in S$, $Q$ has the property 
    \begin{align}
    \label{eq:constructabilityCondition}
      Q_{ii}\ge \sum_{k\in S\backslash\{i\}}Q_{ik}\:\:-\sum_{\{a,b\}\in\binom{S\backslash\{i\}}{2}}Q_{ab}. \qedhere 
    \end{align}
\end{prop}

\begin{proof}
We provide an inductive proof in $|S|$ that is a straightforward application of the principle of Inclusion and Exclusion. Let $F_1,\dots,F_q$ be the facets. For $|S|=2$ we have $Q_{ii}=|F_i\cap F_i|>|F_i\cap F_k|=Q_{ik}$ for any $i$ and $k$. For $|S|=3$,
\begin{align}
Q_{ii}&\ge|\left(F_j\cup F_k\right)\cap F_i|\CR
&=|\left(F_j\cap F_i\right)\cup \left(F_k\cap F_i\right)|\CR
&=|\left(F_j\cap F_i\right)|+|\left(F_k\cap F_i\right)|-|\left(F_j\cap F_k\cap F_i\right)|\CR
&\ge Q_{ji}+Q_{ki}-Q_{jk}
\end{align}
In the last line we have used $|F_j\cap F_k \cap F_i|\le |F_j\cap F_k|=Q_{jk}$.

For the induction step, assume Eqn.(\ref{eq:constructabilityCondition}) is true for $|S|=n$. So for any $i\in S$,
\begin{align}
    Q_{ii}&\ge |\left(\bigcup_{k\in S\backslash\{i\}}F_k\right)\cap F_i|\CR
    &\ge \sum_{k\in S\backslash\{i\}}Q_{ik}\:\:-\sum_{\{a,b\}\in\binom{S\backslash\{i\}}{2}}Q_{ab}
\end{align}

Now for $|S|=n+1$, let $i,j\in S$ be two arbitrary elements, 
\begin{align}
    Q_{ii}&\ge |\left[\left(\bigcup_{k\in S\backslash\{i,j\}}F_k\right)\cup F_j\right]\cap F_i|\CR
    &=|\left[\left(\bigcup_{k\in S\backslash\{i,j\}}F_k\right)\cap F_i\right]\cup \left(F_j\cap F_i\right)|\CR
    &=|\left[\left(\bigcup_{k\in S\backslash\{i,j\}}F_k\right)\cap F_i\right]|+|F_j\cap F_i|-|\left(\bigcup_{k\in S\backslash\{i,j\}}F_k\right)\cap F_j\cap F_i|\CR
    &\ge \left[\sum_{k\in S\backslash\{i,j\}}Q_{ik}\:\:-\sum_{\{a,b\}\in\binom{S\backslash\{i,j\}}{2}}Q_{ab}\right] + Q_{ij} - |\left(\bigcup_{k\in S\backslash\{i,j\}}F_k\right)\cap F_j\cap F_i|\
\end{align}
where we used the induction assumption on $|S\backslash\{j\}|=n$ to get the first term. Further, for the third term,  
\begin{align}
    |\left(\bigcup_{k\in S\backslash\{i,j\}}F_k\right)\cap F_j\cap F_i| & \le |\left(\bigcup_{k\in S\backslash\{i,j\}}F_k\right)\cap F_j|\CR
    &\le \sum_{k\in S\backslash\{i,j\}}Q_{kj}
\end{align}
Therefore, 
\begin{align}
    Q_{ii}   &\ge \left[\sum_{k\in S\backslash\{i,j\}}Q_{ik}+Q_{ij}\right]  \:\:-\sum_{\{a,b\}\in\binom{S\backslash\{i,j\}}{2}}Q_{ab} - \sum_{k\in S\backslash\{i,j\}}Q_{kj}\CR
    &=\sum_{k\in S\backslash\{i\}}Q_{ik}-\sum_{\{a,b\}\in\binom{S\backslash\{i\}}{2}}Q_{ab}. \qedhere
\end{align}
\end{proof}

\section{Clique Complexes}
\subsection{Characterizing Clique Complexes}
\label{sec:cliquecomplexes}
All clique complexes are simplicial complexes, and all simplicial complexes are homeomorphic to clique complexes through subdivision \cite{lickorish1999simplicial, rourke2012introduction}. However, not all simplicial complexes are isomorphic to clique complexes. For example, $K = \{\{1,2\},\{1,3\},\{2,3\}\}$ is not a clique complex because the clique $\{1,2,3\}$ is not in the complex.   In this section we prove a necessary and sufficient condition for a simplicial complex to be a clique complex. In the example for the complex $K = \{\{1,2\},\{1,3\},\{2,3\}\}$, if we take the union of the pairwise intersections of the three facets of $K$, we find 
\begin{align}
    (\{1,2\}\cap\{1,3\})\cup(\{1,2\}\cap\{2,3\})\cup(\{1,3\}\cap\{2,3\}) = \{1,2,3\}.
\end{align}
But the resulting set, $\{1,2,3\}$ is not a face of $K$. Imposing the criteria that the union of the pairwise intersection of any 3-subset of the facets of the complex be a face in the complex turns out to give a necessary and sufficient condition for a simplicial complex to be a clique complex. 

\begin{thm}
\label{thm:cliquecondition}
Let $K$ be a simplicial complex with $n$ facets $\{F_1,\dots,F_n\}$. Then $K$ is a clique complex if and only if $F_{ijk}\equiv(F_i\cap F_j) \cup (F_i\cap F_k) \cup (F_j\cap F_k)$ is a face in $K$ for all $\{i,j,k\}\subseteq [n]$.
\end{thm}
\begin{proof}
In the forward direction, we will show that if $K$ is a clique complex, then every $F_{ijk}$ as defined above is a face of $K$. Let $G$ be the underlying graph to whom $K$ is the clique complex. 
By assumption $F_1, \dots, F_n$ are the maximal cliques of $G$. Take any $\{i,j,k\}\subset [n]$, if any of the intersections $F_i\cap F_j, F_i\cap F_k, F_j\cap F_k$ is empty, then $F_{ijk}$ is also empty (and so a face), so we assume that for a given $i,j,k$  the three pairwise intersections are non-empty. Every vertex of the set $F_i\cap F_j$ is connected to every vertex in $F_i\cap F_k$ because all of these vertices are in the facet $F_i$. Similarly, every vertex in $F_i\cap F_j$ is connected to every element of $F_j\cap F_k$, and every vertex in $F_i\cap F_k$ is connected to every vertex in $F_j\cap F_k$. Therefore, every vertex in $F_{ijk}$ is connected to every other vertex, meaning $F_{ijk}$ is a clique in the graph $G$. If $K$ is a clique complex of $G$, it has to contain all cliques of $G$ as faces, therefore $F_{ijk}$ is a face. 

In the reverse direction, recall that a complex is a clique complex (or flag complex) if and only if its minimal non-faces have size $2$. We will show that if every $F_{ijk}$ is a face then the complex is a flag complex. Suppose not, i.e., there existed a minimal non-face $F$ with size $d \geq 3$. Without loss of generality, let $F = \{1, 2, \ldots, d\}$. Then $f_1 = (F-\{1\})$, $f_2 = (F-\{2\})$, $f_3 = (F-\{3\})$ are all faces. Let $F_1, F_2, F_3$ be the facets that contain $f_1, f_2, f_3$ respectively. By our assumption, $F_{123}$ is a face. And since $f_{123} \subseteq F_{123}$, $f_{123}$ is also a face. But $f_{123}=(F-\{1,2\})\cup(F-\{1,3\})\cup(F-\{2,3\}) = F$ which is a contradiction. Thus, we have shown that no such $F$ exist and indeed $K$ is a flag complex.
\end{proof}

\begin{cor}
\label{cor:cliqueconditionUnion}
Let $K$ be a simplicial complex with $n$ facets $\{F_1,\dots,F_n\}$. Then $K$ is a clique complex if and only if $\bar{F}_{ijk}\equiv(F_i\cup F_j) \cap (F_i\cup F_k) \cap (F_j\cup F_k)$ is either the empty set or a face in $K$ for all $\{i,j,k\}\subseteq \{1,2,\dots,n\}$.
\end{cor}

\begin{proof}
For any three sets $A,B,C$ using de Morgan's law 
\begin{align}
    (A\cap B)\cup(A\cap C)\cup(B\cap C) &=\left[\left(\mathcal{U}-(A \cup B) \right)\cup \left(\mathcal{U}-(A \cup C) \right) \cup \left(\mathcal{U}-(B \cup C) \right)\right]^c\nonumber\\
    &=\left(\mathcal{U}-(A \cup B) \right)^c\cap \left(\mathcal{U}-(A \cup C) \right)^c\cap \left(\mathcal{U}-(B \cup C) \right)^c\CR
    &=(A\cup B)\cap(A\cup C)\cap (B\cup C),
\end{align}
where $\mathcal{U}$ is the universal set and $A^c$ is the complement of set $A$. Therefore, $\bar{F}_{ijk}$ can replace $F_{ijk}$ everywhere in the proof of Thm. \ref{thm:cliquecondition}.
\end{proof} 

\begin{cor}
\label{cor:generalizedCliqueCond}
Let $K$ be a clique complex with $n$ facets $\{F_1,\dots,F_n\}$ and let $\{i_1,\dots,i_\ell\}\subseteq[n]$. Define the set $F_{i_1\dots i_\ell}$ as
    \begin{align}
        F_{i_1\dots i_\ell}&\equiv \bigcup_{\{j_1,\dots,j_{\ell-1}\}\subset\{i_1,\dots,i_\ell\}}\left(\bigcap_{k\in\{j_1,\dots,j_{\ell-1}\}}F_k \right)
    \end{align}
    Then this set is a face in the complex $K$.
\end{cor}
\begin{proof} The statement follows from taking every 3-subset of $\{i_1,\dots,i_\ell\}$, imposing the condition of Theorem (\ref{thm:cliquecondition})  and taking the conjunction. 
\end{proof}
\begin{exmp}
Take a set of four facets $F_1,\dots,F_4$, then the corollary implies that 
\begin{align}
    F_{1234} &= \left(F_1\cap F_2\cap F_3\right)\cup\left(F_1\cap F_2\cap F_4\right)\cup\left(F_1\cap F_3\cap F_4\right)\cup\left(F_2\cap F_3\cap F_4\right)
\end{align}
is a face in the complex.

\end{exmp}

\subsection{Representability of Pure Clique Complexes}

We focus here on pure clique complexes, those whose facets have equal cardinality. The additional structure of 2-determinability of clique complexes makes them such that with an exception to be discussed below, they can be represented uniquely by their facet-adjacency matrix up to isomorphisms. 

One way to study the structure of simplicial complexes is by mapping them to the intersection pattern of some other combinatorial objects. In fact, simplicial complexes that can be represented through the intersection pattern of convex sets as \emph{nerve complexes} have interesting structures studied in the literature (see \cite{tancer2013intersection} for a review). The study of clique complexes as intersection patterns of matrioids is also done, for example here \cite{KASHIWABARA20071910}, where it is shown that any clique complexes with $n$ vertices can  be represented as the intersection of $n-1$ matroids. In this subsection the intersection pattern that we consider is merely the intersection sizes of the facets of a complex. We want to know the minimum extent of knowledge of such intersection sizes needed to uniquely construct a pure complex up to isomorphism. 

Define the \emph{intersection data of degree $k$} of a facet-labeled complex $K$ with $q$ facets and $1 \leq k \leq q $ to be the (labeled) values of $|F_1 \cap \cdots \cap F_k|$ for all choices of $k$ distinct facets $F_1, \ldots, F_k$. Define the \emph{intersection data of $K$} to mean the intersection data over all $k$ in the range $1 \leq k \leq q$. With this framework, we can see that the facet-adjacency matrix contains precisely the intersection data of degrees $1$ and $2$.

\begin{lem}
    \label{lem:intersection-data} A facet-labeled simplicial complex $K$ is uniquely determined (up to isomorphism) from the intersection data of $K$. 
\end{lem}

\begin{proof} The intersection data of $K$ are all the quantities of the form 
\[
|F_{i_1} \cap \cdots F_{i_j}|.
\]
Using this data and inclusion-exclusion, it is easy to compute all $2^q$ quantities of the form $$|(F_1 \text{ or } \overline{F_1})  \cap \cdots \cap  (F_q \text{ or } \overline{F_q})|.$$
This is precisely the vertex data of the complex, which allows us to reconstruct the complex up to the relabeling of the vertices.
\end{proof}

The full intersection data is sufficient to determine a complex, but is it necessary? In light of Theorem~\ref{thm:cliquecondition}, it is reasonable to conjecture that the intersection data up to degree $3$ is sufficient to uniquely reconstruct a pure clique complex. We show that this is not true; in fact, we can say more:

\begin{prop}
\label{prop:constructability-impossible} For any fixed $k$, there exists two non-isomorphic pure clique complexes with the same intersection data up to degree $k$.
\end{prop}

\begin{proof}
The idea of the proof is that it is possible to create a pair of non-isomorphic simplicial complexes (without the pure clique complex conditions) from facets such that their intersection data agree up to degree $k$ but disagree for degree $(k+1)$. We then ``patch'' such a construction to create a pair of non-isomorphic pure clique complexes. 

We will make two facet-labeled simplicial complexes $K$ and $K'$ such that both $K$ and $K'$ have $(k+1)$ labeled facets. $K$ (resp. $K'$) will have facets $F_1,\ldots, F_{k+1}$ (resp. $F_1', \ldots, F_{k+1}'$). We do so by specifying their vertex data. Furthermore, we impose the constraint that these integers equal for any two subsets of the same size (in vertex data notation, $c(a_1\cdots a_q) = c(b_1\cdots b_q)$ whenever the sequences $\{a_1, \ldots, a_q\}$ and $\{b_1, \ldots, b_q\}$ contain the same number of $1$'s). Formally, this means we can dictate that the vertex data (and thus the intersection data) be determined by $k+1$ numbers $f_1, \ldots, f_{k+1}$ such that for all $1 \leq j \leq k+1$ and any $(a_1, \ldots, a_{k+1})$ where exactly $j$ of the $a_i$'s equal $1$, we have $c(a_1\cdots a_q) = f_j.$ The setup also holds for $K'$ with some $k+1$ numbers $f_1', \ldots, f_{k+1}'$ and each $f_j$ replaced by $f_j'$ and $F_i$ replaced by $F_i'$. 

With these very symmetric assumptions, we re-frame our goals as the following conditions:
\begin{enumerate}
\item For all sets of indices $S = \{i_1, \ldots, i_j\}$ for some $j$, we define \[s_S = |F_{i_1} \cap \cdots \cap F_{i_j}|,\]
and similarly $s_S'$ for $K'$.
\item For all $j \leq k$ and $S$ such that $|S| = j$, we have $s_S = s_S'$.
\item For at least one $S$ with $|S| = (k+1)$, we have $s_S \neq s_S'$ (right now there is only one such $S$, but we will end up enlarging $K$ and $K'$).
\end{enumerate}

For now, because of the symmetry of the $f_i$, we see that $s_S$ only depends on the cardinality of $S$. That is, there exist constants $s_1, \ldots, s_{k+1}, s_1', \ldots, s_{k+1}'$ such that if $|S| = j$, $s_S = s_j$. Moreover, we can express each $s_j$ as a linear expression in terms of the $f_i$. By elementary set theory,
\[
s_j = \sum_{i=0}^{k-j} {k-j\choose i} f_{k-i}.
\]
As an example with $k=3, k+1=4$, we get the equations
\begin{align*}
|F_1 \cap F_2 \cap F_3 \cap F_4| & =  s_4 =  f_4 \\
|F_1 \cap F_2 \cap F_3| = \cdots & =  s_3 =  f_3 + f_4 \\
|F_1 \cap F_2| = \cdots & = s_2 = f_2 + 2f_3 + f_4 \\
|F_1| = \cdots & = s_1 =  f_1 + 3 f_2 + 3f_3 + f_4 
\end{align*}
and analogous equations for $K'$. 

There are many solutions that satisfy our constraints, but we use the following especially interesting one: set $f_i = 1$ for odd $i$ (and $f_i = 0$ otherwise) and set $f_i' = 1$ for even $i$ (and $f_i' = 0$ otherwise). It is easy to check then $s_i = s_i'$ is satisfied for all $i$ except $s_{k+1} \neq s'_{k+1}$. 

Finally, we must also check the conditions of Proposition~\ref{prop:vertex-data-reconstruction}. For our specific construction, we notice:
\begin{itemize}
    \item Since $f_1 = 1$ and $k+1 \geq 3$, this implies that for any $i \neq j$, choosing $a_k = 0$ for all $k$ except $a_i = 1$ gives $c(a_1\cdots a_{k+1}) = 1$, satisfying Proposition~\ref{prop:vertex-data-reconstruction}.
    \item Since $f'_2 = 1$ and $k+1 \geq 3$, this implies that for any $i \neq j$, we can choose $a_i = 1$, $a_j = 0$, and $a_k = 0$ for all $k$ except one arbitrary index. We also have $c(a_1 \cdots a_{k+1}) = 1$ for this sequence.
\end{itemize}
This argument ensures that our construction of the intersection data corresponds to actual complexes.

Now, suppose we come up with any such pair of complexes $K$ and $K'$. Our goal is to turn them into pure clique complexes. First, let $v$ be the bigger number of vertices between $K$ and $K'$. Now, we add facets (both of size $v$) $F_{k+2}$ to $K$ and $F_{k+2}'$ to $K'$ such that they contain all vertices in their corresponding complexes, padding with new vertices that only belong to the new facet when needed. We can check that:
\begin{itemize}
    \item for any $S \subset \{1, \ldots, k+1\}$ with $|S| \leq k$, our desired conditions remain unchanged when we compare the intersection data indexed by $S$ for $K$ and $K'$;
    \item for any $S$ such that $(k+2) \in S$ and $|S| \leq k$, we have
    \[
    |F_{i_1} \cap \cdots \cap F_{i_j = k+2}| = |F_{i_1} \cap \cdots \cap F_{i_{j-1}}| = |F_{i_1}' \cap \cdots \cap F_{i_{j-1}}'| = |F'_{i_1} \cap \cdots \cap F_{i_j=k+2}'|.
    \]
\end{itemize}
In other words, our conditions remain met, though it is no longer true that $s_S$ only depends on $|S|$; for each $|S|$ there are now two possible values of $s_S$, one for when $S$ includes $(k+2)$ and one for when it does not.
   
Our pair of complexes now both have $(k+2)$ facets, though the facets have different numbers of vertices. For the final step, in $K$, add $v-s_1$ new vertices to each facet $F_1, \ldots, F_{k+1}$ only belonging to that facet, and perform the analogous operation on $K'$ (recall that $s_1 = s_1'$). Now every facet has $v$ vertices and we have added $(k+1)(v-s_1)$ new vertices to each complex, we can check that:
\begin{itemize}
    \item the intersection data of degree $1$ changed for both $K$ and $K'$, but they changed all the $|F_i| = |F_i'|$ to $v$ (it helps to recall that $v-s_1 = v-s_1'$ because $s_1 = s_1'$).
    \item the intersection data of any degree $k > 1$ remains unchanged (since none of the new vertices are involved in any pairwise intersections or higher).
\end{itemize}

Finally, it remains to prove that $K$ is a pure clique complex. It suffices to show that any clique in the $1$-skeleton for $K$ is a subset of one of the facets $F_1, \ldots, F_{k+2}$. To see this, first note that no such clique can involve the $(k+1)(v-s_1)$ ``new'' vertices we just added (because they only belong to a single facet), so the clique must belong entirely to the first $v$ vertices. But then by construction they belong to $F_{k+2}$, so we are done (as the analogous result holds for $K'$).
\end{proof}

\begin{exmp}
We give an explicit example of our construction for $k=3$. Solving our equations, we obtain the facet-incidence matrices
\[
B = \begin{bmatrix}
1 & 1 & 1 & 0 & 1 & 0 & 0 & 0 \\
1 & 1 & 0 & 1 & 0 & 1 & 0 & 0 \\
1 & 0 & 1 & 1 & 0 & 0 & 1 & 0 \\
0 & 1 & 1 & 1 & 0 & 0 & 0 & 1
\end{bmatrix}, B' = \begin{bmatrix}
1 & 1 & 1 & 0 & 0 & 0 & 1 \\
1 & 0 & 0 & 1 & 1 & 0 & 1 \\
0 & 1 & 0 & 1 & 0 & 1 & 1 \\
0 & 0 & 1 & 0 & 1 & 1 & 1
\end{bmatrix}.
\]

These correspond to simplicial complexes (which we call $K$ and $K'$) involving $4$ facets of size $4$; call these facets $F_1, \ldots, F_4$ and $F_1', \ldots, F_4'$. It is easy to check that their intersection data matches up to degree $3$ (for example, $|F_1 \cap F_2| = |F_1' \cap F_2'| = 2$, as do any pairwise intersections of facets) but differ at degree $4$ ($|F_1 \cap \cdots \cap F_4| = 0$, but $|F_1' \cap \cdots \cap F_4'| = 1$).

Because $\max(8,7) = 8$, we add a facet to each complex (call them $F_5$ and $F_5'$ respectively) that uses all the existing vertices and has size $v = 8$. For $F_5'$, this means adding one more vertex that is not in any of $F_1'$ through $F_4'$.

Finally, for each of the $8$ facets ($F_1, \ldots, F_4, F_1', \ldots, F_4'$), add $4$ new vertices only belonging to that facet. This brings all facets to the same size. Concretely, we now have $2$ facet-incidence matrices
\[
\begin{bmatrix}
1 & 1 & 1 & 0 & 1 & 0 & 0 & 0 & 1 & 1 & 1 & 1 & 0 & 0 & 0 & 0 & 0 & 0 & 0 & 0 & 0 & 0 & 0 & 0 \\
1 & 1 & 0 & 1 & 0 & 1 & 0 & 0 & 0 & 0 & 0 & 0 & 1 & 1 & 1 & 1 & 0 & 0 & 0 & 0 & 0 & 0 & 0 & 0 \\
1 & 0 & 1 & 1 & 0 & 0 & 1 & 0 & 0 & 0 & 0 & 0 & 0 & 0 & 0 & 0 & 1 & 1 & 1 & 1 & 0 & 0 & 0 & 0 \\
0 & 1 & 1 & 1 & 0 & 0 & 0 & 1 & 0 & 0 & 0 & 0 & 0 & 0 & 0 & 0 & 0 & 0 & 0 & 0 & 1 & 1 & 1 & 1 \\
1 & 1 & 1 & 1 & 1 & 1 & 1 & 1 & 0 & 0 & 0 & 0 & 0 & 0 & 0 & 0 & 0 & 0 & 0 & 0 & 0 & 0 & 0 & 0
\end{bmatrix}
\]
and
\[
\begin{bmatrix}
1 & 1 & 1 & 0 & 0 & 0 & 1 & 0 & 1 & 1 & 1 & 1 & 0 & 0 & 0 & 0 & 0 & 0 & 0 & 0 & 0 & 0 & 0 & 0 \\
1 & 0 & 0 & 1 & 1 & 0 & 1 & 0 & 0 & 0 & 0 & 0 & 1 & 1 & 1 & 1 & 0 & 0 & 0 & 0 & 0 & 0 & 0 & 0 \\
0 & 1 & 0 & 1 & 0 & 1 & 1 & 0 & 0 & 0 & 0 & 0 & 0 & 0 & 0 & 0 & 1 & 1 & 1 & 1 & 0 & 0 & 0 & 0 \\
0 & 0 & 1 & 0 & 1 & 1 & 1 & 0 & 0 & 0 & 0 & 0 & 0 & 0 & 0 & 0 & 0 & 0 & 0 & 0 & 1 & 1 & 1 & 1 \\
1 & 1 & 1 & 1 & 1 & 1 & 1 & 1 & 0 & 0 & 0 & 0 & 0 & 0 & 0 & 0 & 0 & 0 & 0 & 0 & 0 & 0 & 0 & 0
\end{bmatrix}.
\]
We leave it to the reader to check that these are indeed facet-incidence matrices for $8$-pure simplicial complexes, that the intersection data up to degree $3$ are the same, but the intersection data for degree $4$ differ (specifically, we still have $|F_1 \cap \cdots \cap F_4| = 0$ and $|F_1' \cap \cdots \cap F_4'| = 1$).
\end{exmp}

Given these results, it is interesting to ask, ``are there classes of complexes for whom it is possible to get the full intersection data from the facet-adjacency matrix?'' In other words, are there classes of pure clique complexes that can be uniquely represented by their facet-adjacency matrix (up to isomorphisms)? Consider any $3$ facets, $F_1,F_2,F_3$. We say that they have \emph{triangle intersection} if the intersections $F_{\hat{1}}\equiv\bar{F_1}\cap F_2\cap F_3,\;F_{\hat{2}}\equiv F_1\cap \bar{F_2}\cap F_3,\; F_{\hat{3}}\equiv F_1\cap F_2\cap \bar{F_3}$ are all non-empty. We will prove that any pure clique complex which has no triangle intersection is uniquely determined from its facet-adjacency matrix (in other words, its intersection data of degrees $1$ and $2$). 

\begin{thm}
    \label{thm:possible} 
Suppose $K$ is $p$-pure with $q$ facets and we know that $K$ has no set of $3$ facets with triangle intersection. Then the intersection data of degrees $1$ and $2$ determines $K$ uniquely up to isomorphism. 
\end{thm}

\begin{proof}
We prove by induction that we can construct the intersection data of degree $n$ for all $n \leq q$. As base cases, we already have $n = 1$ and $n=2$ data in the facet-adjacency matrix.

We now show how, assuming the intersection data for up to $n$ gives the $(n+1)$-st layer as well, as long as $n+1 \leq q$. Consider any set of $(n+1)$ facets $F_1, \ldots, F_{n+1}$. Then we know the following $(n+1)$ pieces of degree-$n$ intersection data
\[
f_{\hat{i}} = |F_1 \cap \cdots \cap  \hat{F_i} \cap \cdots \cap F_{n+1}| 
\]
as $i$ runs from $1$ to $(n+1)$, 
where $\hat{F_i}$ denotes that we omit $F_i$. Now, we claim that
\begin{align}
\label{eq:minimizationRule}
|F_1 \cap \cdots \cap F_{n+1}| = \min(f_{\hat{1}}, \ldots, f_{\hat{n+1}}).
\end{align}
This is automatically true if any of the $f_{\hat{i}}$ equal 0 so we will consider the case when $f_{\hat{i}}>0$ for all $i$. Suppose the opposite, i.e., $|F_1 \cap \cdots \cap F_{n+1}| < \min(f_{\hat{1}}, \ldots, f_{\hat{n+1}})$, then since $(F_1 \cap \cdots \cap F_{n+1})$ is a subset of the set corresponding to every $f_{\hat{i}}$, we know 
\[
|F_1 \cap \cdots \cap F_{n+1}| < f_{\hat{i}}
\]
for all $i$. Then, since for all $i$
\[ |F_1 \cap \cdots \cap F_{n+1}| + |F_1 \cap \cdots \overline{F_i} \cdots \cap F_{n+1}| = |F_1 \cap \cdots \hat{F_i} \cdots \cap F_{n+1}| = f_{\hat{i}},
\]
we must have $|F_1 \cap \cdots \overline{F_i} \cdots \cap F_{n+1}| \geq 1$ for all $1\le i\le n+1$. We now know that there must exist $(n+1)$ vertices $v_1, \ldots, v_{n+1}$ of $K$ such that each $v_i$ belongs to all $F_1, \ldots, F_{n+1}$ and does \textbf{not} belong to $F_i$. Take any three of these vertices, WLOG, $v_1,v_2,v_3$. Then, since $\{v_1\}\in \bar{F_1}\cap F_2\cap \dots\cap F_{(n+1)} \in \bar{F_1}\cap F_2\cap F_3$, $\{v_2\}\in F_1\cap \bar{F_2}\cap\dots\cap F_{(n+1)} \in F_1\cap \bar{F_2}\cap\ F_3$, and $\{v_3\}\in F_1\cap F_2\cap \bar{F_3}\cap\dots F_{(n+1)}\in F_1\cap F_2\cap \bar{F_3}$, the three facets $F_1, F_2, F_3$ form a triangle intersection, a contradiction. 

Rewinding, we have concluded that $F_1 \cap \cdots \cap F_{n+1}$ is uniquely determined by intersection data up to $n$, for any facets $F_1, \ldots, F_{n+1}$. In other words, the intersection data of degree $n+1$ are deterministic given the intersection data of degree $n$, as long as the clique complex has no triangle intersection. This finishes the inductive step and thus the induction.
\end{proof}

\begin{exmp}
For an example of what goes wrong when we remove the requirement of the absence of triangle intersection, consider the following two 4-pure complexes: 
\begin{align}
    K_1&=\{\{1, 2, 4, 5\},\{1, 3, 6, 7\},\{2, 3, 8, 9\},\{1, 2, 3, 10\}\},\CR
    K_2&=\{\{1, 2, 5, 6\},\{1, 3, 7, 8\},\{1, 4, 9, 10\},\{1, 2, 3, 4\}\}
\end{align}

\begin{figure}
\centering
    \includegraphics[width=0.9\textwidth]{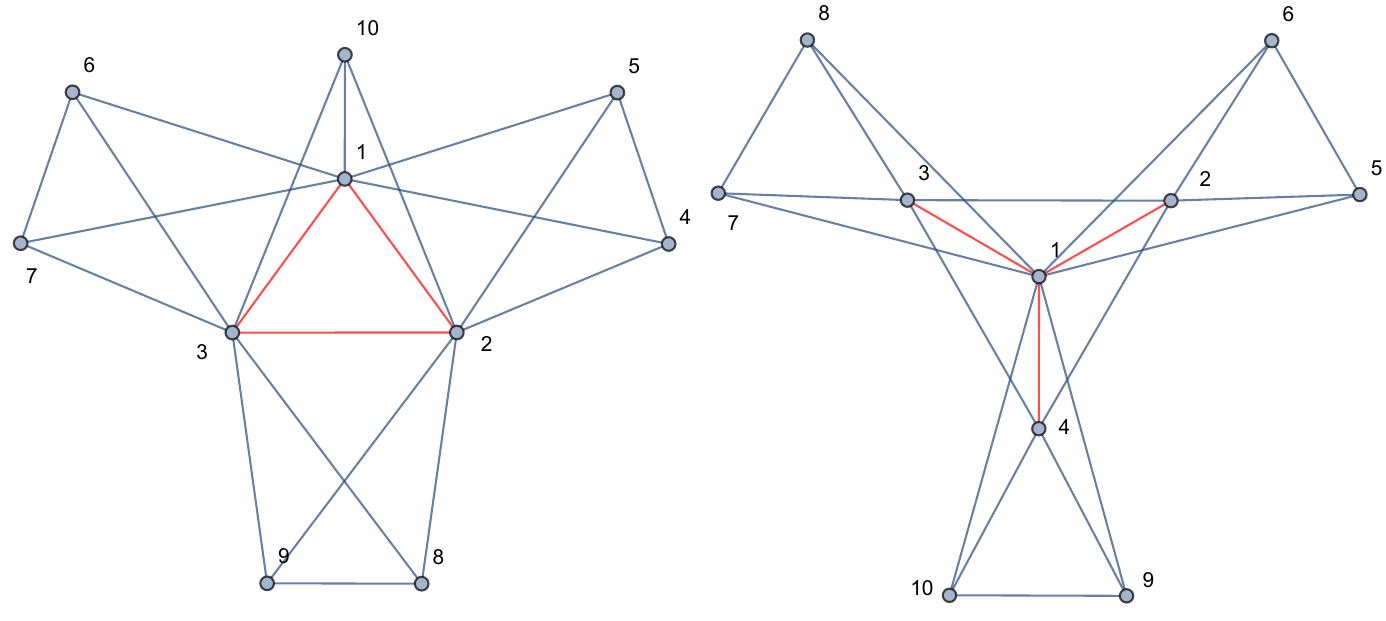} 
\caption{Two non-isomorphic 4-pure clique complexes with the same facet-adjacency matrix. The triangle and claw intersections are highlighted in red.}
\label{fig:tetrahedronfree}
\end{figure}
The two complexes are drawn in Fig.(\ref{fig:tetrahedronfree}). One can check that these they have identical facet-adjacency matrix
\[Q=\left(
\begin{matrix}
4 & 1 & 1 & 2 \\
1 & 4 & 1 & 2 \\
1 & 1 & 4 & 2 \\
2 & 2 & 2 & 4
\end{matrix}\right).
\]
\end{exmp}

\section{Counting Pure Complexes}
\label{sec:counting}
In this section we will  count the number of vertex-labeled pure simplicial complexes and provide bounds on the number of vertex-labeled pure clique complexes for a fixed purity $p$ with and facet order $q$.
Let $\mathcal{S}_v(p,q,n)$ be the set of vertex-labeled $p-$pure simplicial complexes with $q$ facets and $n$ vertices so that 
 \begin{align}
     \mathcal{S}_v(p,q) = \bigcup_{n}\mathcal{S}_v(p,q,n)
 \end{align} 
 Let $s_v(p,q,n)$ and $s_v(p,q)$ be the cardinalities of the respective sets:
 \begin{align}
s_v(p,q,n)=|\mathcal{S}_v(p,q,n)|,\;\; s_v(p,q)=|\mathcal{S}_v(p,q)| =\sum_n s_v(p,q,n)
 \end{align} 
Similarly let $\mathcal{W}_v(p,q,n)$ be the set of $p$-pure clique complexes of order $q$ with $n$ labeled vertices so that $\mathcal{W}(p,q)=\cup_n \mathcal{W}(p,q,n)$. Let the cardinalities of the sets be similarly denoted $w_s(p,q,n)$ and $w_s(p,q)$. Table \ref{tab:purecounts} below gives $s_v(p,q)$ and $w_v(p,q)$ for the first few $p$ and $q$.

\begin{table}[H]
    \centering  
\begin{tabular}{|c||c|c|c|c|c|c|}
\hline\
$q$&1&2&3&4&5&6\\
\hhline{|=||=|=|=|=|=|=|}
   $s_v(2,q)$ & 1 & 6 & 62 & 900 & 16,824 &384,668 \\
   \hline
   $w_v(2,q)$ & 1 & 6 & 61 & 878 & 16,323 &371,782\\
   \hline\hline
   $s_v(3,q)$ & 1 & 31 & 2,649 &441,061&121,105,865 &49,615,422,851 \\
  \hline
   $w_v(3,q)$ & 1 & 31 & 2,495 &394,920&104,268,613 &41,419,848,444\\
     \hline\hline
     $ s_v(4,q) $& 1& 160 & 116,360 & 231,173,330 & 974,787,170,226 &7,500,396,185,804,060 \\
     \hline
     $ w_v(4,q) $& 1&160 &101,875 & 178,682,745 & 679,213,720,913 & 4,793,115,687,225,971\\
     \hline
\end{tabular}
\caption{The number of vertex-labeled pure simplicial $s_v(p,q)$ and clique $w_v(p,q)$ complexes for $p\in\{2,3,4\}$ and $q\in \{1,\dots,6\}$}
\label{tab:purecounts}
\end{table}
At the time of the completion of this paper, except for the first row none of the sequences  in the table above appear in the Online Encyclopedia of Integer Sequences (OEIS). The first row of numbers (i.e., $s_v(2,q)$) is also the number of labeled graphs with no isolated vertices and $q$ edges studied in \cite{bender1997asymptotic} and appearing in OEIS \cite{oeisA121251}. The authors of \cite{bender1997asymptotic} gave only an assymptotic formula for these sequences. The OEIS entry \cite{oeisA121251} gives the a formula for $s_v(2,q)$ as an infinite sum
\begin{align}
    s_v(2,q)=\sum_{m\ge 0}\frac{1}{2^{(m+1)}}\binom{\binom{m}{2}}{q}
\end{align}
A quick numerical check reveals that the infinite sum formula generalizes to higher purity values as 
\begin{align}
    s_v(p,q)=\sum_{m\ge 0}\frac{1}{2^{(m+1)}}\binom{\binom{m}{p}}{q}
\end{align}
In the section below we give a formula for the number $s_v(p,q,n)$ as a simple finite sum and $s_v(p,q)$ as a double sum. 

\subsection{Counting Pure Simplicial Complexes}
There is a simple expression for the number of vertex-labeled pure simplicial complexes $s_v(p,q,n)$ provided in the proposition below.

\begin{prop}
The number of $p$-pure simplicial complexes of facet-order $q$ with $n$ vertices $s_v(p,q,n)$ is given by
\begin{align}
\label{eq:svpqn}
    s(p,q,n)&=\sum_{k=0}^{n}(-1)^{k+n}\binom{n}{k}\binom{\binom{k}{p}}{q}
\end{align}
\end{prop}

\begin{proof}
$s_v(p,q,n)$ is the cardinality of the set of $q$-combinations of the set of all $p$-combinations of $[n]$ that cover $[n]$. The maximum number of vertices that a pure simplicial complex in $\mathcal{S}_v(p,q)$ can have is $pq$ corresponding to the case where all $q$ facets are disconnected. On the other hand, the minimum number of vertices that any $X\in \mathcal{S}_v(p,q)$ can have, $n_0$, corresponds to the edge-maximal graph with $q$ facets. That will be the smallest complete graph having $q$ or more $p$-cliques. Since the complete graph $K_n$ over $n$ vertices has $\binom{n}{p}$ $p-$cliques, $n_0$ is given by 
 \begin{align}
     n_0=\min\{n | \binom{n}{p}\ge q\}
 \end{align} 

In general if we consider the set of all $q$-combinations of the $p$-combinations of $[n]$, this set can be split according to the number of vertices $k$ covered by the union of facets. There are $\binom{n}{k}$ ways to choose $k$ vertices out of $[n]$ which are covered by the union of the facets. For example, for $p=2,q=2, n=3$, a vertex labeled 2-pure complexes with 2 facets and 3 vertices is an element of 2-combinations of $\binom{[3]}{2}=\{\{1,2\},\{1,3\},\{2,3\}\}$. Taking all 2-combinations of $\binom{[3]}{2}$ gives $3=\binom{3}{2}$ 2-pure complexes of facet order 3. Therefore, \begin{align}
    \binom{\binom{n}{p}}{q}=\sum_{k=0}^{n}\binom{n}{k}s_v(p,q,k)
\end{align}
Then, the Binomial inversion formula  gives the desired expression for $s_v(p,q,n)$.

\end{proof}

\begin{cor}
    The number of $p$-pure vertex-labeled  simplicial complexes with $n$ vertices is given by 
    \begin{align}
        s_v(p,n)=\sum_{k=0}^n(-1)^{n+k}\binom{n}{k}2^{\binom{k}{p}}
    \end{align}
\end{cor}
\begin{proof}
    \begin{align}
        s_v(p,n)&=\sum_{q=q_{min}}^{q_{max}}s_v(p,q,n)\CR
        &=\sum_{q=q_{min}}^{q_{maz}}\sum_{k=0}^{n}(-1)^{n+k}\binom{n}{k}\binom{\binom{k}{p}}{q}\CR
        &=\sum_{k=0}^n\sum_{q=0}^{\binom{k}{p}}(-1)^{n+k}\binom{n}{k}\binom{\binom{k}{p}}{q}\CR
        &=\sum_{k=0}^n (-1)^{n+k}\binom{n}{k} 2^{\binom{k}{p}}
    \end{align}
\end{proof}
\begin{rem}
    Again utilizing the binomial inversion formula we find \begin{align} 2^{\binom{n}{p}}=\sum_{k=0}^{n}\binom{n}{k}s_v(p,k)
    \end{align}
    This is saying that summing over all $k$ of the product of the number of ways of choosing $k$ labels from $n$ possibilities and the number of vertex labeled $p$-pure complexes with $k$ vertices gives the total number of subsets of $\binom{[n]}{p}$. Specifically when $p=2$, the number of all simple graphs with $n$ vertices is $2^{\binom{n}{2}}$. That equals the sum over all $k$ of $\binom{n}{k}$ times the number of all simple graphs with $k$ non-isolated vertices. So we are counting the number of all graphs with $n$ vertices by summing over $k$ of the number graphs with exactly $(n-k)$ isolated vertices. 
\end{rem}

\subsection{Connection to the number of alignments}
There is an interesting combinatorial connection between the number of pure simplicial complexes and the enumeration of \emph{alignments} as discussed in \cite{slowinski1998number}. We define an \emph{alignment} of $k$ strips of length $p$ to be a $k \times m$ matrix where:
\begin{enumerate}
    \item Each element of the matrix has value $0$ or $1$;
    \item there are $p$ $1$'s in every row;
    \item no column has all $0$'s.
\end{enumerate}

An example of an alignment of $3$ strips of length $4$ is
\[
\begin{bmatrix}
    1 & 1 & 0 & 0 & 1 & 1 \\
    1 & 1 & 1 & 1 & 0 & 0 \\
    0 & 1 & 0 & 1 & 1 & 1
\end{bmatrix}.
\]
Let $f(p,k)$ be the number of alignments of $k$ strips of length $p$. The work \cite{slowinski1998number} is motivated by counting the number of ways of aligning $k$ DNA strands of length $p$ alongside each other where the strands are bound in particular places to each other. 
\begin{cor}
    The number of $p$-pure vertex-labeled  simplicial complexes with $q$ facets is given by 
    \begin{align}
        s_v(p,q)=\frac{1}{q!}\sum_{k=0}^q s(q,k)f(p,k),
    \end{align}
    where $s(q,k)$ are the Stirling numbers of the first kind and $f(p,k)$ is the number of alignments of $k$ strips of length $p$. 
\end{cor}

\begin{proof}

Algebraically, we can compute 
    \begin{align}
        s_v(p,q)&=\sum_{n=0}^{pq}s_v(p,q,n)\CR
        &=\sum_{n=0}^{pq}\sum_{k=0}^n(-1)^{n+k}\binom{n}{k}\binom{\binom{k}{p}}{q}\CR
        &=\frac{1}{q!}\sum_{n=0}^{pq}\sum_{k=0}^n(-1)^{n+k}\binom{n}{k}\left(\binom{k}{p}\right)_{q}\CR
        &=\frac{1}{q!}\sum_{n=0}^{pq}\sum_{k=0}^n(-1)^{n+k}\binom{n}{k}\sum_{j=0}^{q}s(q,j)\binom{k}{p}^j\CR
        &=\frac{1}{q!}\sum_{j=0}^qs(q,j)\sum_{n=p}^{pj}\sum_{k=0}^n(-1)^{n+k}\binom{n}{k}\binom{k}{p}^j\CR
        &=\frac{1}{q!}\sum_{j=0}^qs(q,j)f(p,j),
    \end{align}
    where 
    \begin{align}
    f(p,j)=\sum_{N=p}^{pj}\sum_{i=0}^{N}(-1)^i\binom{N}{i}\binom{N-i}{p}^j 
    \end{align}
    is the number of alignments of $j$ strips of length $p$ as found in \cite{slowinski1998number}.

We now also provide a combinatorial interpretation of this equality. First, define $t(p,q) = q! s_v(p,q)$. Notice that $t$ counts the number of vertex-labeled $p$-pure clique complexes with $q$ different \textbf{labeled} facets. We may rewrite our result as
\[
 t(p,q) = \sum_{j=0}^q s(q,j) f(p,j).
\]
By construction of Stirling numbers, this is equivalent to
\[
f(p, q) = \sum_{j=0}^q S(q,j) t(p, j),
\]
where $S(q,j)$ are the Stirling numbers of the second kind. It suffices to show that the left-hand-side and the right-hand-side counts the same thing. The following steps outline a bijection between these two objects:

\begin{enumerate}
\item First, we look at the left-hand-side. This corresponds to alignments of $q$ strips of length $p$. 
\item Now, for any row vector $v$ that appears at rows $r_1, \ldots, r_s$, we represent these rows by a single ``representative row" $v$ labeled by the set $\{r_1, \ldots, r_s\}$ 
\item Thus, alignments of $q$ strips of length $p$ are in bijection with sets of \textbf{unique} ``representative rows" whose labels form a partition of $[q]$.
\item Given a set of unique ``representative rows" $\{s_1, \ldots, s_k\}$ labeled by $S_1, \ldots, S_k$ respectively, we can construct a unique $p$-pure complex in the following way: for each strip $s_i$, pick the $p$ column indices that appear in the row vector of $s_i$, and construct a $p$-facet with them.
\item Finally, note that the right-hand-side exactly counts the number of $p$-pure clique complexes whose facets are labeled by a partition of $[q]$: we first pick how many facets, $j$, to have, then assign the labels in $S(q,j)$ ways.
\end{enumerate}

As an example of this bijection, take the following alignment counted by $f(3,4)$:

\[
\begin{bmatrix}
    1 & 1 & 1 & 0 & 0 & 0 \\
    0 & 0 & 1 & 1 & 1 & 0 \\
    0 & 0 & 1 & 1 & 1 & 0 \\
    1 & 0 & 0 & 0 & 1 & 1 \\
\end{bmatrix}.
\]

This alignment can be represented as $3$ labeled ``representative strips'':
\begin{itemize}
    \item $(1, 1, 1, 0, 0, 0)$ labeled by $\{1\}$;
    \item $(0, 0, 1, 1, 1, 0)$ labeled by $\{2,3\}$;
    \item $(1, 0, 0, 0, 1, 1)$ labeled by $\{4\}.$
\end{itemize}

This set of labeled representative strips then corresponds to the following $3$-pure complex with labeled facets:
\begin{itemize}
    \item $\{1,2,3\}$ (labeled by $\{1\}$);
    \item $\{3,4,5\}$ (labeled by $\{2,3\}$);
    \item $\{1,5,6\}$ (labeled by $\{4\}$). \qedhere
\end{itemize}
\end{proof}

\subsection{Upper bound on the number of Pure Clique Complexes}
A modification of Eqn.(\ref{eq:svpqn}) provides an upper bound on $w_v(p,q,n)$. If $W\in \mathcal{W}_v(p,q,n)$ is a clique complex, then the minimum possible value of the number of its vertices,  $n^*$, is achieved when the graph formed is edge maximal $p-$pure graph of clique order $q$. Edge-maximal $p-$pure graphs are Tur\'an graphs \cite{chao1982maximally}. A $p$-pure Tur\'an graph with $m$ vertices is a $p-$partite complete graph $K_{q+1,q+1,\dots,q,q}$ where the $m$ vertices are distributed as equitably as possible among the $p$ partitions so that $m = q p+r$, with $r<p$. The number $q$ of maximal $p-$cliques in a Tur\'an graph with $m$ vertices is given by  
\begin{align}
\label{eqn:Turan}
    q=\prod_{i=1}^p \left\lfloor\frac{m+i-1}{p}\right\rfloor
\end{align} 
Given $p$ and $q$, the minimum number of nodes $n^*$ is given by the minimum $m$ that gives $q$ from equation (\ref{eqn:Turan}) above. Let us define the {\em Tur\'an number} $r(p,q)$ as follows:

\begin{defn}
The Tur\'an number $r(p,q)$ of a $p-$pure complex with $q$ facets is the minimum number of vertices needed to construct a $p-$pure Tur\'an graph with $q$ facets.
\begin{align}
    r(p,q)=\min \{n|\prod_{i=1}^p \left\lfloor\frac{n+i-1}{p}\right\rfloor\ge q\}
\end{align}
\end{defn}
Table (\ref{tab:turanNumbers}) gives Tur\'an numbers $r(p,q)$ for $p=2,3,4,5$ and a range of clique orders $q$.
\begin{table}[h]
    \centering
    \begin{tabular}{|c||c|c|c|c|}
    \hline
        $k$&$r(2,q)$&$r(3,q)$&$r(4,q)$&$r(5,q)$\\
        \hline\hline
        1 & 2 &3 &4 & 5 \\
        \hline
        2 & 3&4&5&6\\
        \hline
        3&4&5&6&7\\
        \hline
        4&4&5&6&7\\
        \hline
        5&5&6&7&8\\
        \hline
        6&5&6&7&8\\
        \hline
        7&6&6&7&8\\
        \hline
        8&6&6&7&8\\
        \hline
        9&6&7&8&9\\
        \hline
        10&7&7&8&9\\
        \hline
    \end{tabular}
    \caption{Tur\'an number $r(p,q)$ for $p=2,3,4,5$ and $q=1,\dots,10$.}
    \label{tab:turanNumbers}
\end{table}

Therefore, summing $s_v(p,q,n)$ over $n\in \{r(p,q),r(p,q)+1,\dots,pq\}$ gives an upper bound on $w_v(p,q)$,
\begin{align}
    w_v(p,q)\le\sum_{n=r(p,q)}^{pq} s_v(p,q,n)
\label{eq:numcliquecomplexes}
\end{align}

\section{Summary of Results}
\label{sec:conclusion}
We derived a number of results regarding simplicial and clique complexes with a fixed number of facets. Any vertex-labeled simplicial complex with $q$facets can be represented uniquely through a list of its facets, through the facet-incidence matrix $B$, and through the vertex data map $c:\{0,1\}^q\rightarrow\{1,2,\dots\}$. Given the facet-incidence matrix $B$ of a simplex, its facet-adjacency matrix $Q=BB^T$ is a $q\times q$ symmetric integer matrix. In section 2 we considered the converse case and established criteria to determine when a given $q\times n$ matrix $B$ of 1's and 0's is realizable as the facet-incidence matrix of a simplicial complex. We also characterized maps $c$ and square symmetric integer matrix $Q$ that can be realized as vertex data and facet-adjacency of a simplicial complex. 

Focusing attention to clique complexes, we characterized clique complexes purely in terms of the list of their facets through Thm.(\ref{thm:cliquecondition}) which established a necessary and sufficient condition on the list of facets of a simplicial complex to be a clique complex. This provided a means to establish two results on the unique representability of pure clique complexes from their intersection data. We proved that it is in general not possible to represent all pure clique complexes with level $k$ intersection data for any fixed $k<q$. However, any pure clique complex that is triangle-intersection free can be uniquely represented from its level 2 intersection data up to isomorphism. 

In the last section we counted the number of $p$-pure simplicial complexes of facet-order $q$ with $n$ labeled vertices, 
\begin{align}
    s_v(p,q,n)&=\sum_{k=0}^{n}(-1)^{k+n}\binom{n}{k}\binom{\binom{k}{p}}{q}
\end{align} 
We established a connection between the number of $p$-pure simplicial complexes with labeled vertices and facets $t(p,q)=q!s_v(p,q)$ and the number of alignments $f(p,q)$ of $q$ strips of length $p$ that commonly appear in enumerations of the number of DNA alignments, 
\begin{align}
    f(p,q)&=\sum_{j=0}^qS(q,j)t(p,j)
\end{align}
Finally, in Eqn.(\ref{eq:numcliquecomplexes}) we placed an upper bound on the number of $p$-pure clique complexes of facet-order $q$.

\section*{Acknowledgements}
This work is supported in part by the National Science Foundation LEAPS-MPS grant (award number 2138323) and DOE RENEW-HEP grant (award number DE-SC0024518). CE is also supported by the McNair fellowship program at San Jose State University. We want to thank Boris Alexeev, Kevin Iga, and Allen Liu for valuable feedback on drafts of the paper.
\pagebreak
\appendix

\printbibliography 

\end{document}